\documentclass[12pt,a4paper]{amsart}
\usepackage{a4wide}
\usepackage[english]{babel}
\usepackage[T2A]{fontenc}
\usepackage[utf8]{inputenc} 
\usepackage{amsfonts}
\usepackage{amssymb, amsthm, amscd}
\usepackage{amsmath}
\usepackage{mathtools}
\usepackage{needspace}
\usepackage{etoolbox}
\usepackage{lipsum}
\usepackage{comment}
\usepackage{cmap}
\usepackage[pdftex]{graphicx}
\usepackage[unicode]{hyperref}
\usepackage[matrix,arrow,curve]{xy}
\usepackage[usenames,dvipsnames]{xcolor}
\usepackage{colortbl}
\usepackage{textcomp}
\usepackage{cite}
\usepackage{euscript}

\newcommand{\Z}{\mathbb{Z}}

\newcommand{\N}{\mathbb{N}}

\newcommand{\M}{\mathfrak{M}}
\newcommand{\A}{\mathfrak{A}}
\newcommand{\B}{\mathfrak{B}}

\newcommand{\eps}{\varepsilon}

\newcommand{\sub}{\subseteq}

\renewcommand{\ge}{\geqslant}
\renewcommand{\le}{\leqslant}
\newcommand{\sm}{\setminus}
\newcommand{\map}[3]{#1\colon #2\to #3}
\newcommand{\phan}{\vphantom{,}}

\renewcommand{\>}{\rangle}

\newcommand{\emp}{\varnothing}
\newcommand{\SR}{\mathrm{SR}}
\newcommand{\ASR}{\mathrm{ASR}}

\newcommand{\AFSR}{\mathrm{AFSR}}

\renewcommand{\t}{\reflectbox{t}}

\newcommand{\SL}{\mathop{\mathrm{SL}}\nolimits}
\newcommand{\GL}{\mathop{\mathrm{GL}}\nolimits}

\newcommand{\tc}{\text{,}}
\newcommand{\tp}{\text{.}}
\renewcommand{\tilde}{\widetilde}

\DeclareMathOperator{\Max}{Max}
\DeclareMathOperator{\IUm}{IUm}
\DeclareMathOperator{\BSdim}{BSdim}
\DeclareMathOperator{\Jdim}{Jdim}

\DeclareMathOperator{\hgt}{ht}
\DeclareMathOperator{\Ann}{Ann}
\DeclareMathOperator{\res}{res}
\DeclarePairedDelimiter{\ceil}{\lceil}{\rceil}
\DeclarePairedDelimiter{\floor}{\lfloor}{\rfloor}

\theoremstyle{plain}
\newtheorem{thm}{Theorem}[section]
\newtheorem{lem}[thm]{Lemma}
\newtheorem{prop}[thm]{Proposition}

\theoremstyle{definition}
\newtheorem{defn}[thm]{Definition}
\newtheorem{cor}[thm]{Corollary}

\theoremstyle{remark}
\newtheorem*{rem}{Remark}

\AtBeginEnvironment{thm}{\begin{samepage}}
	\AtEndEnvironment{thm}{\end{samepage}}
\AtBeginEnvironment{lem}{\begin{samepage}}
	\AtEndEnvironment{lem}{\end{samepage}}
\AtBeginEnvironment{st}{\begin{samepage}}
	\AtEndEnvironment{st}{\end{samepage}}
\AtBeginEnvironment{crit}{\begin{samepage}}
	\AtEndEnvironment{crit}{\end{samepage}}
\AtBeginEnvironment{ax}{\begin{samepage}}
	\AtEndEnvironment{ax}{\end{samepage}}
\AtBeginEnvironment{defn}{\begin{samepage}}
	\AtEndEnvironment{defn}{\end{samepage}}
\AtBeginEnvironment{cor}{\begin{samepage}}
	\AtEndEnvironment{cor}{\end{samepage}}
\AtBeginEnvironment{note}{\begin{samepage}}
	\AtEndEnvironment{note}{\end{samepage}}
\AtBeginEnvironment{prop}{\begin{samepage}}
	\AtEndEnvironment{prop}{\end{samepage}}

\makeatletter
\def\@settitle{\begin{center}%
		\baselineskip14\p@\relax
		\bfseries
		\@title
	\end{center}%
}

\def\@evenhead{\hfil\sc Pavel Gvozdevsky\hfil}
\def\@oddhead{\hfil\sc Bounded reduction\hfil}
\makeatother

\title{Bounded reduction of orthogonal matrices over polynomial rings}

\keywords{Split orthogonal group, surjective stability of $K_1$, bounded reduction, polynomial rings}
\subjclass[2020]{19B14(primary), 19B10, and 20G35(secondary)} 

\author{Pavel Gvozdevsky}
\date{}
\address{Chebyshev Laboratory, St. Petersburg State University, 14th Line V.O., 29B, Saint Petersburg 199178 Russia}

\begin{document}

\maketitle

\begin{abstract}
	We prove that a matrix from the split orthogonal group over a polynomial ring with coefficients in a small-dimensional ring can be reduced to a smaller matrix by a bounded number of elementary orthogonal transformations. The bound is given explicitly. This result is an effective version of the early stabilisation of the orthogonal $K_1$ functor proven by Suslin and Kopeiko in \cite{SusKopLGPSO}. Since the similar effective results for special linear and symplectic groups are obtained by Vaserstein in \cite{VasBounded}, the present paper closes the problem for split classical groups.
\end{abstract}

\section{Introduction}
This paper deals with the last remaining case of bounded reduction in split classical groups over polynomial rings. Namely we study the case of orthogonal groups.

Chevalley groups over certain rings have bounded width with respect to the elementary generators. For example this holds for Dedekind domains of arithmetic type, see \cite{CarterKeller}, \cite{CarterKellerZ}, \cite{Morris}, \cite{TavgenChevalley}, \cite{TavgenTwisted},\cite{MorganRapinchukSury},\cite{ErovenkoRapinchuk2001},\cite{ErovenkoRapinchuk2006}. Results on such {\it bounded generation} are of great value, for example they are connected to the {\it congruence subgroup property}, see \cite{Lubotzky},\cite{PlatonovRapinchuk}, and to Margulis--Zimmer conjecture, see \cite{ShalomWillis}. However bounded generation occurs very rarely in the sense that classes of rings for which it is known to hold are pretty narrow. Nevertheless, for some applications it is enough to have a weaker result, such as: bounded length of conjugates of elementary generators (see \cite{StepVavDecomp}), bounded length of commutators (see \cite{SivStep},\cite{StepVavLength},\cite{HSVZwidth}), or bounded generation with respect to a larger set of generators. Bounded reduction is a variation of the last property. 

A given Chevalley group $G$ over a given ring is said to have bounded reduction if any element of $G$ can be decomposed as a product of bounded number of elementary generators and one (not necessarily elementary) element from a certain subsystem subgroup. In other words, it means that one can reduce any element to the subsystem subgroup by bounded number of elementary transformations. Without requirement for the number of elementary transformations to be bounded this property is called the surjective stability of the $K_1$-functor. In papers \cite{SteinStability}, \cite{SteinMats}, \cite{PlotProc}, \cite{PlotkinStability}, \cite{PlotkinE7}, and \cite{GvozK1Surj} this problem is considered for rings that satisfy certain conditions on stable rank, absolute stable rank, or other similar conditions. Actually, from the proofs of the theorems in these papers one can recover the bound on the required number of elementary transformations, despite the fact that this bound is not stated in papers explicitly. Therefore, these are results on bounded reduction.

However, conditions on stable rank are still very strong. Even though small Jacobson dimension implies small stable rank, rings with large Jacobson dimension usually fail to have small stable rank. In the present paper, we consider another important class of rings. Namely we take a polynomial ring in arbitrary number of variables with coefficients in a small-dimensional ring. Here we use Krull dimension because the techniques require for dimension to behave well with respect to adding an independent variable. 

Without the bound on the number of elementary transformations similar result for classical groups is known. This is so called early surjective stability of the $K_1$-functor. For the special linear group this was proved by Suslin in \cite{SuslinSerreConj}. Similar result for the orthogonal group follows from \cite{SusKopLGPSO}, and for the symplectic group it is proven in \cite{KopejkoSpStab}, see also \cite{GMVLGPSp}, \cite{KopejkoSpDimOne}. Note that if the ring of coefficient is a Dedekind domain or a smooth algebra over a field, then this result for all Chevalley groups follows from the homotopy invariance of the non-stable $K_1$-functor, see \cite{AbeLGP},\cite{StavPolynom},\cite{StavSerreConj}.

In the case of special linear and symplectic groups, there are similar results for Laurent polynomial rings, see \cite{KopejkoLorSL}, \cite{KopejkoLorSp}.

In the paper \cite{VasBounded}, Vaserstein obtained the effective version of the Suslin result, i.e. he proved the bounded reduction for the special linear group over a polynomial ring, and gave this bound explicitly. From this result he deduced that the elementary subgroup of the general linear group over an arbitrary finitely generated commutative ring has Kazhdan's property (T).

In \cite{Shalom}, the basic connection between bounded generation and property (T) has been established and used to estimate the Kazhdan Constants for $\SL_n(\Z)$. Later the bounds for these constants were improved in \cite{KasabovKazdanSL}. In order to deduce property (T) from the Vaserstein's result one needs to refer to \cite{SholomAlgebraisation}. 

In fact, property (T) for Chevalley groups and groups similar to them has already been studied by other methods, see \cite{PropertyTforUniversalLattices},\cite{EJKPropertyT}. However, we believe that the bounded reduction has an independent value, and we aim to study this question for other Chevalley groups. It was noted in the concluding remarks of \cite{VasBounded} that the bounded reduction for the symplectic group follows formally from the case of special linear group. Therefore, as we said, the last remaining classical case is the one of orthogonal groups. The main result of the present paper is the following theorem.

{\thm\label{main} Let $C$ be a commutative Noetherian ring and $\dim C=D<\infty$. Let $A=C[x_1,\ldots,x_n]$. Then for any $r\ge \max(3,D+2)$, every matrix from the split orthogonal group $O(2r,A)$ resp. $O(2r+1,A)$ can be reduced to the subgroup
$$
\begin{pmatrix}
	1 & & \\
	 &O(2r-2,A) & \\
	 &        & 1
\end{pmatrix}\qquad\text{resp.}\qquad\begin{pmatrix}
1 & & \\
&O(2r-1,A) & \\
&        & 1
\end{pmatrix}
$$
by multiplication from the left by 
\begin{gather*}
	N=n(11r-7)+\left(nD+\tfrac{r(r-1)}{2}\right)\left(\tfrac{r(r-1)}{2}+\ceil{\tfrac{r-1}{2}}+8r-2\right)+10r-10\\
	\text{resp.}\\  N=n(12r-8)+\left(nD+\tfrac{r(r-1)}{2}\right)\left(\tfrac{r(r-1)}{2}+\ceil{\tfrac{r-1}{2}}+9r-2\right)+11r-9
\end{gather*}
elementary orthogonal transvections.}

\smallskip

Therefore, on one hand this theorem is an effective version of the result from \cite{SusKopLGPSO}, and on the other hand, an extension of \cite{VasBounded} to the orthogonal group.  

Very soon we reduce Theorem \ref{main} to Theorem \ref{ColumnReduction}, which deals with an arbitrary column vector rather than a column of a particular matrix. So we make a convention that the word "column" means an arbitrary column vector.

In fact, Theorem \ref{ColumnReduction} below is something stronger than just bonded reduction. It is easy to see that we can deduce from it the following corollary, which is a special case of Theorem 7.10 of \cite{SusKopLGPSO}.

{\cor Let $C$ be a commutative Noetherian ring and $\dim C=D<\infty$, $A=C[x_1,\ldots,x_n]$. Let $V$ be an $A$-module equipped with a quadratic form $q$. Assume that an orthogonal sum of $(V,q)$ with a hyperbolic plane is a free module with a split quadratic form of rank $2r$ or $2r+1$ with $r\ge \max(3,D+2)$. Then $V$ is a free module and $q$ is split.}

\smallskip

The paper is organised as follows. In Section \ref{PreliminariesAndNotation}, we give all necessary preliminaries and introduce basic notation. In Section \ref{AbsoluteFlexibleStableRank}, we introduce the new notion of an absolute flexible stable rank, and prove some of its properties. In Sections \ref{ReductionToPropositions}, \ref{MonicSection}, and \ref{EliminationSection} we give the proof of the main result.

\section{Preliminaries and notation}
\label{PreliminariesAndNotation}
\subsection*{Rings, ideals and dimensions}

By a ring we always mean associative and commutative ring with unity. 

\smallskip

If $R$ is a ring, then by $R^*$ we denote the set of invertible elements in $R$.

\smallskip

For the elements $r_1$,$\ldots$,$r_k\in R$, we denote by $\<r_1,\ldots,r_k\>$ the ideal in $R$ generated by these elements.
 
\smallskip

In the present paper we use three different notions of a ring dimension.

\begin{itemize}

\item By $\dim R$ we denote Krull dimension of the ring $R$. That is the supremum of the lengths of all chains of prime ideals.

\smallskip

\item By $\Jdim R=\dim\Max(R)$ we denote the dimension of the maximal spectrum $\Max(R)$ of the ring $R$. It is equal to the supremum of the lengths of all chains of such prime ideals that coincide with its Jacobson radical.

\smallskip

\item By $\BSdim R$ we denote the Bass--Serre dimension of a ring $R$. That is the minimal $\delta$ such that $\Max(R)$ is a finite union of irreducible Noetherian subspaces of dimension not greater than $\delta$.

\end{itemize}

Obviously, for a Noetherian ring $R$ we have
$$
\BSdim R\le \Jdim R \le \dim R\tp
$$

\smallskip

The following property of Bass--Serre dimension is well known; see Lemma 4.17 in \cite{BakNonabelian}.

{\lem\label{BSInduction} Let $R$ be a ring with $\BSdim R=\delta<\infty$. Then it has a finite collection $P_1$,$\ldots$,$P_m$ of maximal ideals, such that for any element $s\in R\sm \bigcup_i P_i$ we have $\BSdim R/(s)<\delta$. In case where $\delta=0$, this means that $s\in R^*$.}

\subsection*{Even split orthogonal group}

Let $R$ be a commutative ring and $V$ be a free $R$-module of rank $2r$ with basis $e_1$,$\ldots$,$e_r$,$e_{-r}$,$\ldots$,$e_{-1}$; this ordering of the basis elements will be used throughout the paper. Consider the quadratic form $q$ on the module $V$ defined as follows
$$
q(x)=\sum_{i=1}^r x_ix_{-i}\tc
$$
where $x_i$ are the coordinates of $x$ in our basis.

The group of automorphism of $V$ preserving the form $q$ is called the {\it even split orthogonal group} and denoted by $O(2r,R)$. We identify its elements with their matrices in the basis $e_1$,$\ldots$,$e_r$,$e_{-r}$,$\ldots$,$e_{-1}$.

For an element $\xi\in R$ and indices $i,j=1,\ldots,-1$ with $i\ne\pm j$, we set
$$
T_{i,j}(\xi)=T_{-j,-i}(-\xi)=e+\xi e_{i,j}-\xi e_{-j,-i}\tp
$$
Such matrices are called {\it elementary orthogonal transvections}. It is easy to see that $T_{i,j}(\xi)\in O(2r,R)$.

The subgroup of $O(2r,R)$, generated by all elementary orthogonal transvections is called the {\it elementary subgroup} and denoted by $EO(2r,R)$.

For any $k\in\N$, we denote by $EO(2r,R)^{\le N}$ the subset of $EO(2r,R)$ consisting of matrices that can be expressed as the product of no more than $N$ elementary orthogonal transvections.

\subsection*{Odd split orthogonal group}

Now let $V$ be a free $R$-module of rank $2r+1$ with basis $e_1$,$\ldots$,$e_r$,$e_0$,$e_{-r}$,$\ldots$,$e_{-1}$, and the quadratic form $q$ on the module $V$ defined as follows
$$
q(x)=x_0^2+\sum_{i=1}^r x_ix_{-i}\tp
$$
The group of automorphism of $V$ preserving the form $q$ is called the {\it odd split orthogonal group} and denoted by $O(2r+1,R)$.

In the odd case, we have two types of elementary orthogonal transvections. First type is the same as in the even case:
$$
T_{i,j}(\xi)=T_{-j,-i}(-\xi)=e+\xi e_{i,j}-\xi e_{-j,-i}\qquad i,j=1,\ldots,-1\tc \; i\ne\pm j \tc
$$
and we refer to them as {\it long root transvections}. Second type is called {\it short root transvections}, and those are the matrices
$$
T_{i,0}(\xi)=e+2\xi e_{i,0}-\xi e_{0,-i}-\xi^2 e_{i,-i}\qquad i=\pm 1,\ldots\pm r\tp
$$

{\rem The terms "long root" and "short root" come from the terminology for Chevalley groups, which are defined by root systems; see, for example, \cite{Humphreys}.

The even split orthogonal group corresponds to the root system
$$
D_n=\{e_i\pm e_j\colon i\ne j, 1\le i,j\le n\}\tc
$$
in which all the roots have the same length.

However, the odd split orthogonal group corresponds to the root system
$$
B_n=\{e_i\pm e_j\colon i\ne j, 1\le i,j\le n\}\cup\{e_i\colon 1\le i\le n\}\tc
$$
in which there are long roots (of length $\sqrt{2}$) and short roots (of length 1).}

\smallskip

As in the even case, we denote by $EO(2r+1,R)$ the subgroup of $O(2r+1,R)$, generated~by all elementary orthogonal transvections, and by  $EO(2r+1,R)^{\le N}$, the subset of $EO(2r+1,R)$ consisting of matrices that can be expressed as product of no more than $N$ elementary orthogonal transvections.

Obviously, the group $O(2r,R)$ embeds into $O(2r+1,R)$ as the stabiliser of $e_0$.

\subsection*{Matrix transpose with respect to antidiagonal}

Since we use weird numeration of rows and columns in orthogonal matrices, we need a notation for the matrix transpose with respect to antidiagonal. Such transpose of a matrix $g$ we denote by $g^{\t}$.

\smallskip

By $g^t$ we denote regular transpose of the matrix $g$.

\subsection*{Hyperbolic embedding}

Consider the map
$$
\map{H}{\GL(r,R)}{O(2r,R)}\qquad g\mapsto \begin{pmatrix}
g & \\  & (g^{-1})^{\t}	
\end{pmatrix}\tc
$$

This map is called the {\it hyperbolic embedding}.

\smallskip

Note that $H$ takes elementary transvections in $\GL(r,R)$ to elementary orthogonal transvections. Here by elementary transvections in $\GL(r,R)$ we mean the matrices
$$
t_{i,j}(\xi)=e+\xi e_{i,j}\tc
$$

where $\xi\in R$, $i,j=1,\ldots,r$, $i\ne j$.

In the odd case, we denote by $H$ the composition
$$
\xymatrix{
\GL(r,R)\ar[r]^H & O(2r,R) \ar@{^{(}->}[r] & O(2r+1,R)\tp
}
$$

\subsection*{ASR-condition}
Recall that a commutative ring $R$ satisfies the absolute stable rank condition $\ASR_d$ if for any row $(b_1,\ldots,b_d)$ with coordinates in $R$, there exist elements $c_1$,$\ldots$,$c_{d-1}\in R$ such that every maximal ideal of $R$ containing the ideal $\<b_1+c_1b_d,\ldots,b_{d-1}+c_{d-1}b_d\>$ contains already the ideal $\<b_1,\ldots,b_d\>$. This notion was introduced in \cite{EstesOhm} and used in \cite{SteinStability}, \cite{SteinMats} and then in \cite{PlotProc},\cite{PlotkinStability}, \cite{PlotkinE7}, and \cite{GvozK1Surj} to study stability problems.

If we assume that a row $(b_1,\ldots,b_d)$ is unimodular, then the absolute stable rank condition boils down to the usual stable rank condition $\SR_d$ (see \cite{Bass64},\cite{VasersteinStability}).

Absolute stable rank satisfies the usual properties, namely for every ideal $I\unlhd R$ condition $\ASR_d$ for $R$ implies $\ASR_d$ for the quotient $R/I$, and if $d\ge d'$, then $\ASR_{d'}$ implies $\ASR_d$. Finally, it is well known that if the maximal spectrum of $R$ is a Noetherian space of dimension $\Jdim R=d-2$, then both conditions $\ASR_d$ and $\SR_d$ are satisfied (see \cite{EstesOhm},\cite{MagWandarKalVas},\cite{SteinStability}).

\section{Absolute flexible stable rank}
\label{AbsoluteFlexibleStableRank}
In this section, we introduce a new type of stable rank condition called the {\it absolute flexible stable rank}. Here is the definition.

{\defn A commutative ring $A$ satisfies the absolute flexible stable rank condition $\AFSR_d$ if for any row $(b_1,\ldots,b_d)$ with coordinates in $A$, there exists an element $c_1\in A$ such that for any invertible element $\eps_1\in A^*$, there exists $c_2\in A$ such that for any $\eps_2\in A^*$, $\ldots$, there exists $c_{d-1}\in A$ such that for any $\eps_{d-1}\in A^*$, every maximal ideal of $A$ containing the ideal $\<b_1+\eps_1c_1b_d,\ldots,b_{d-1}+\eps_{d-1}c_{d-1}b_d\>$ contains already the ideal $\<b_1,\ldots,b_d\>$. }

\smallskip

One can think of it as follows. Two players are playing a game. Player 1 chooses a row $(b_1,\ldots,b_d)$ with coordinates in $A$. Then they take turns starting with Player 2. Player 2 in his $i$-th turn chooses element $c_i\in A$; after that Player 1 in his turn chooses invertible element $\eps_i\in A^*$. Player 2 wins if after $d$ turns every maximal ideal of $A$ containing the ideal $\<b_1+\eps_1c_1b_d,\ldots,b_{d-1}+\eps_{d-1}c_{d-1}b_d\>$ contains already the ideal $\<b_1,\ldots,b_d\>$. A commutative ring $A$ satisfies the absolute flexible stable rank condition $\AFSR_d$ if Player 2 has a winning strategy.

Now we show that the condition $\AFSR_d$ holds for small-dimensional rings. That generalises the result of \cite{EstesOhm}.

{\lem\label{DimmaxImpliesAFSR} Let $A$ be a commutative ring. Assume that $\Max(A)$ is Noetherian and $\Jdim A\le d-2$. Then $A$ satisfies $\AFSR_d$.}
\begin{proof}
	Following \cite{EstesOhm}, for any ideal $I\unlhd A$ we denote by $J(I)$ the intersection of all the maximal ideals containing $I$, and we set $J=\{I\unlhd A\mid J(I)=I\}$. The ideals from $J$ are in correspondence with the closed subsets of $\Max(A)$.
	
	For an ideal $I\unlhd A$, a prime ideal $P\in J$ is called a {\it component} of $I$ if $P$ is minimal among the primes of $J$ that contains $I$. The assumption for $\Max(A)$ to be Noetherian implies that any ideal has only finitely many components.
	
	Now we show that making his $i$-th turn Player 2 can guarantee that any component of ideal $I=\<b_1+\eps_1c_1b_d,\ldots,b_{i-1}+\eps_{i-1}c_{i-1}b_d\>$ (for $i=1$ set $I=0$) that contains $b_i+\eps_ic_ib_d$ also contains $b_d$, and that doing so is a winning strategy.
	
	When the $i$-th turn of Player 2 starts, the ideal $I$ is already fixed and has finitely many component. Let $X$ be the set of components of $I$ that does not contain $b_d$. All Player 2 has to do is to guarantee that $b_i+\eps_ic_ib_d$ is not contained in the union of ideals from $X$.
	
	Set $Y=\{P\in X\mid b_i\in P\}$. Since the components of $I$ cannot contain each other, we can choose for any $\A$,$\B\in X$ an element $z(\A,\B)\in \A\setminus \B$. Then to achive his goal Player 2 can set
	$$
	c_i=\sum_{\B\in Y}\prod_{\A\in X\setminus\{\B\}} z(\A,\B)\tp 
	$$
	
	Let us prove that Player 2 won. Let $b_i'=b_i+\eps_ic_ib_d$, $i=1,\ldots,d-1$. Assume that some maximal ideal $\M$ contains $\<b_1',\ldots,b_{d-1}'\>$, but does not contain $b_d$. We claim that in this case there exists a chain of primes in $J$
	$$
	P_0\subsetneq P_1\subsetneq\ldots\subsetneq P_{d-1}=\M\tc
	$$
	which contradicts the assumption on $\Jdim A$. 
	
	We build this chain from the end in such a way that $\<b_1',\ldots,b_i'\>\le P_i$ but $b_{i+1}'\notin P_i$. We set $P_{d-1}=\M$. Then for $i<d-1$ we choose $P_i$ to be the component of $\<b_1',\ldots,b_i'\>$ that is contained in $P_{i+1}$, such a component exists by Zorn's lemma. Since $P_i\sub \M$, it follows that $b_d\notin P_i$, and by construction of $b_{i+1}'$, this implies that $b_{i+1}'\notin P_i$.
	
	The fact that $b_{i+1}'\notin P_i$ guarantees that the inclusions in the chain are proper.
\end{proof}

The following lemma shows how one can use the $\AFSR$ condition.

{\lem\label{UseAFSR} Let $A$ be a commutative ring, and $S$ be a multiplicative system in $A$. Assume that the localisation $A[S^{-1}]$ satisfies $\AFSR_d$. Then for any row $(b_1,\ldots,b_d)$ with coordinates in $A[S^{-1}]$ and for any $s\in S$, there exist $c_1$,$\ldots$,$c_{d-1}\in sA$ such that such that every maximal ideal of $A[S^{-1}]$ containing the ideal $\<b_1+c_1b_d,\ldots,b_{d-1}+c_{d-1}b_d\>$ contains already the ideal $\<b_1,\ldots,b_d\>$.}

\begin{proof}
	Just let Player 1 choose his invertible element in such a way that the resulting coefficient be in $sA$.
\end{proof}

\section{Reduction of Theorem \ref{main} to propositions \ref{MakeMonic} and \ref{KillTheVariable}}
\label{ReductionToPropositions}

We denote by $\IUm_{2r}A$ the set isotropic unimodular columns of over a ring $A$:
$$
\IUm_{2r}A=\{b\in A^{2r}\colon b\text{ is unimodular, }q(b)=0\}\tp
$$ 

Similarly, for the odd case we denote
$$
\IUm_{2r+1}A=\{b\in A^{2r+1}\colon (b_1,\ldots,b_r,2b_0,b_{-r},\ldots,b_{-1})\text{ is unimodular, }q(b)=0\}\tp
$$ 

We have\footnote{We need to clarify why in the odd case the first column of an orthogonal matrix remains unimodular after multiplying the middle element by two. That is because the symmetric bilinear form associated to $q$ is $Q(x,y)=q(x+y)-q(x)-q(y)=2x_0y_0+\sum_{i=1}^{r}(x_iy_{-i}+x_{-i}y_i)$. If $b$ and $c$ are the first and the last columns of some matrix in $O(2r+1,A)$, then $b_1c_{-1}+\ldots+b_bc_{-r}+(2b_0)c_0+b_{-r}c_r+\ldots+b_{-1}c_1=Q(b,c)=Q(e_1,e_{-1})=1$. So the row $(b_1,\ldots,b_r,2b_0,b_{-r},\ldots,b_{-1})$ is unimodular.} $O(2r,A)e_1\le \IUm_{2r}A$ and $O(2r+1,A)e_1\le \IUm_{2r+1}A$, where $e_1$ is the first column of the identity matrix.

If the first column of a matrix $\alpha\in O(2r,A)$ resp. $O(2r+1,A)$ coincides with $e_1$, then $\alpha$ automatically has shape
$$
\begin{pmatrix}
	1 & *  		& *\\
	  & \beta	& *\\
	  &			&1
\end{pmatrix}\tc
$$
where $\beta\in O(2r-2,A)$ resp. $O(2r-1,A)$. Hence it can be reduced to 
$$
\begin{pmatrix}
	1 & *  		& *\\
	& \beta	& \\
	&			&1
\end{pmatrix}\tc
$$
which automatically should be 
$$
\begin{pmatrix}
	1 &  		& \\
	& \beta		& \\
	&			&1
\end{pmatrix}\tc
$$
 by multiplication from the left by $2r-2$ resp $2r-1$ elementary orthogonal transvections.

\smallskip

Therefore, theorem \ref{main} follows trivially from the following result.

{\thm\label{ColumnReduction} Under the condition of Theorem \ref{main}, for every column $b\in\IUm_{2r}A$ resp. $\IUm_{2r+1}A$ we have
$$
e_1\in EO(2r,A)^{\le N}b\qquad\text{resp.}\qquad e_1\in EO(2r+1,A)^{\le N}b \tc
$$
where 
\begin{gather*}
N=n(11r-7)+\left(nD+\tfrac{r(r-1)}{2}\right)\left(\tfrac{r(r-1)}{2}+\ceil{\tfrac{r-1}{2}}+8r-2\right)+8r-8\\
 \text{resp.}\\  N=n(12r-8)+\left(nD+\tfrac{r(r-1)}{2}\right)\left(\tfrac{r(r-1)}{2}+\ceil{\tfrac{r-1}{2}}+9r-2\right)+9r-8\tp
\end{gather*}}

Consider the lexicographic order on the monomials in variables $x_1$,$\ldots$,$x_n$. That is the order where $x_1^{k_1}\ldots x_n^{k_n}$ is bigger than $x_1^{l_1}\ldots x_n^{l_n}$ if for some $m$ we have $k_i=l_i$ for $i<m$, and $k_m>l_m$. A polynomial in $A=C[x_1,\ldots,x_n]$ is called lexicographically monic if its leading coefficient in lexicographic order is equal to one.

\smallskip

Further we reduce Theorem \ref{ColumnReduction} to the following two propositions.

{\prop\label{MakeMonic} Under the condition of Theorem \ref{main}, for every column $b\in\IUm_{2r}A$ resp. $\IUm_{2r+1}A$ there exists a column
$$
b'\in EO(2r,A)^{\le 11r-7}b\qquad\text{resp.}\qquad EO(2r+1,A)^{\le 12r-8}b\tc
$$
such that its entry $b'_{-2}$ is lexicographically monic.
}

{\prop\label{KillTheVariable}
Let $B$ be a commutative ring such that $\BSdim B=\delta<\infty$. Let $r\ge 3$, $A=B[y]$, $b\in\IUm_{2r}A$ resp. $\IUm_{2r+1}A$ such that its entry $b_{-2}$ is monic. Then
$$
EO(2r,A)^{\le N}b\cap B^{2r}\ne\emp\quad \text{resp.} \quad EO(2r+1,A)^{\le N}b\cap B^{2r+1}\ne\emp\tc
$$
where $N=\delta\left(\tfrac{r(r-1)}{2}+\ceil{\tfrac{r-1}{2}}+8r-2\right)$ resp. $\delta\left(\tfrac{r(r-1)}{2}+\ceil{\tfrac{r-1}{2}}+9r-2\right)$
}

\smallskip

First we need the following lemma.

{\lem\label{Normalisation} Let $f\in C[x_1,\ldots,x_n]$ be a lexicographically monic polynomial. Then there exists an invertible change of variables
$$
x_1,\ldots,x_n \leftrightarrow y_1,\ldots, y_n\tc
$$
such that $f$ becomes monic in $y_n$.}

\begin{proof}
	Take $K>\deg f$. Set $x_i=y_i+y_n^{K^{n-i}}$, $i=1$,$\ldots$,$n-1$, and $x_n=y_n$.
\end{proof}

\smallskip

Now we deduce Theorem \ref{ColumnReduction} from Propositions~\ref{MakeMonic} and~\ref{KillTheVariable}.

\smallskip

Take $b\in\IUm_{2r}A$ resp. $\IUm_{2r+1}A$. By Proposition~\ref{MakeMonic} there exists a column
$$
b'\in EO(2r,A)^{\le 11r-7}b\qquad\text{resp.}\qquad EO(2r+1,A)^{\le 12r-8}b\tc
$$
such that its entry $b'_{-2}$ is lexicographically monic. Applying Lemma~\ref{Normalisation}, we change variables to $y_1$,$\ldots$,$y_n$ so that  $b'_{-2}$ is now monic in $y_n$. Now we apply Proposition~\ref{KillTheVariable} to $B=C[y_1,\ldots,y_{n-1}]$. Note that $\BSdim B\le \dim B=D+n-1$. Hence we can obtain a column from
\begin{gather*}
EO(2r,A)^{\le N_1'}b'\cap B^{2r}\le EO(2r,A)^{\le N_1}b\cap B^{2r}\\
\text{resp.}\\
EO(2r+1,A)^{\le N_1'}b'\cap B^{2r+1}\le EO(2r+1,A)^{\le N_1}b\cap B^{2r+1}\tc
\end{gather*}
where 
$$
N_1'=(D+n-1)\left(\tfrac{r(r-1)}{2}+\ceil{\tfrac{r-1}{2}}+8r-2\right)\quad \text{resp.}\quad (D+n-1)\left(\tfrac{r(r-1)}{2}+\ceil{\tfrac{r-1}{2}}+9r-2\right)\tc
$$
 and $N_1=N_1'+11r-7$ resp. $N_1'+12r-8$.

Repeating this argument $n$ times we can obtain a column from
$$
EO(2r,A)^{\le N_n}b\cap C^{2r}\qquad\text{resp.}\qquad EO(2r+1,A)^{\le N_n}b\cap C^{2r+1}\tc
$$
where
\begin{gather*}
N_n=n(11r-7)+\left(nD+\tfrac{r(r-1)}{2}\right)\left(\tfrac{r(r-1)}{2}+\ceil{\tfrac{r-1}{2}}+8r-2\right)\\\text{resp.}\\ N_n=n(12r-8)+\left(nD+\tfrac{r(r-1)}{2}\right)\left(\tfrac{r(r-1)}{2}+\ceil{\tfrac{r-1}{2}}+9r-2\right)\tp
\end{gather*}

Since $\dim C\le r-2$, the ring $C$ satisfies the condition $\ASR_r$. Hence the proof of surjective stability of $K_1$ for orthogonal groups given in \cite{SteinStability} (Theorem 2.1), which is essentially a proof of bounded reduction under the $\ASR$-condition, implies\footnote{Stein uses action not on columns but on rows, and in the odd case there is a discrepancy with the present paper about where the factor of two arises. So if our column, say is $(v_1,\ldots,v_r,v_0,v_{-r},\ldots,v_{-1})^t$, one has to apply Stein's proof to the row $(v_1,\ldots,v_r,2v_0,v_{-r},\ldots,v_{-1})$. Note that this row remains unimodular.} that applying to our column $6r-6$ resp. $7r-7$ elementary orthogonal transvections we can make its first entry equal to one. Then by $2r-2$ resp. $2r-1$ elementary orthogonal transvections we can annihilate entries from $2$ to $-2$, and the last $-1$-th entry becomes zero automatically because the column remains isotropic.

Therefore,
$$
e_1\in EO(2r,A)^{\le N}b\qquad\text{resp.}\qquad e_1\in EO(2r+1,A)^{\le N}b\tc
$$
where 
\begin{gather*}
	N=n(11r-7)+\left(nD+\tfrac{r(r-1)}{2}\right)\left(\tfrac{r(r-1)}{2}+\ceil{\tfrac{r-1}{2}}+8r-2\right)+8r-8\\
	\text{resp.}\\  N=n(12r-8)+\left(nD+\tfrac{r(r-1)}{2}\right)\left(\tfrac{r(r-1)}{2}+\ceil{\tfrac{r-1}{2}}+9r-2\right)+9r-8\tp
\end{gather*}

\section{Obtaining a monic polynomial}
\label{MonicSection}

In this section, we give the proof of Proposition~\ref{MakeMonic}. First we need some preparation.

\smallskip

Recall that the height $\hgt_A(P)$ of a prime ideal $P$ in a ring $A$ is the is the largest number $h$ (or $\infty$ if such a number does not exist) such that there exists a chain of different prime ideals
$$
P_0\subsetneq P_1\subsetneq\ldots\subsetneq P_n=P\tp
$$

The height $\hgt_A(I)$ of an arbitrary ideal $I\unlhd A$ is defined by
$$
\hgt_A(I)=\inf\{\hgt_A(P)\colon P\text{ is prime, } I\le P\}\tp
$$

{\lem\label{HeightAndMonicPolynomials} Let $C$ be a Noetherian ring, $A=C[x_1,\ldots,x_n]$, $I\unlhd A$ and $\hgt_A(I)\ge \dim(C)+1$. Then the ideal $I$ contains a lexicographically monic polynomial.}

\begin{proof}
	Let $J$ be an ideal in $B=C[x_2,\ldots,x_n]$ consisting of leading terms of polynomials in $I$. Here we consider $A$ as $B[x_1]$. By Lemma~10.4 of \cite{VasSusSerre} we have
	$$
	\hgt_B(J)\le \hgt_A(I)\le \dim(C)+1\tp
	$$
	
	In case $n=1$ where $B=C$, this implies that $1\in J$. Hence $I$ contains a monic polynomial. That is the base of induction. For the induction step, note that by induction hypothesis $J$ contains a lexicographically monic polynomial. Thus $I$ contains a polynomial $f$ with leading term in $x_1$ being lexicographically monic, but this means exactly that $f$ is lexicographically monic itself.
\end{proof}

{\lem\label{DimOfLocalisation} Let $C$ be a Noetherian ring, $A=C[x_1,\ldots,x_n]$. Let $S$ be a multiplicative system of lexicographically monic polynomials in $A$. Then we have $\dim A[S^{-1}]\le \dim C$.}
\begin{proof}
	Assume the converse. Then there exists a chain 
	$$
	P_0\subsetneq P_1\subsetneq\ldots\subsetneq P_k
	$$
	of prime ideals in $A[S^{-1}]$, where $k>\dim C$. Hence we have
	$$
	P_0\cap A\subsetneq P_1\cap A\subsetneq\ldots\subsetneq P_k\cap A\tp
	$$
	Thus $\hgt(P_k\cap A)\ge k\ge \dim C+1$, and by Lemma~\ref{HeightAndMonicPolynomials} the ideal $P_k\cap A$ contains a lexicographically monic polynomial. Hence $P_k=A[S^{-1}]$. This is a contradiction.

\end{proof}

Now we recall a definition from \cite{VasBounded}.

{\defn\label{mu} Let $A$ be an associative ring with 1, $s$ be a central element of $A$, $r\ge 2$, $v\in A^{r-1}$ (a column over $A$), $u\in\phan^{r-1}A$ (a row over $A$). We define an $r$ by $r$ matrix over $A$ by
$$
\mu(u,s,v)=\begin{pmatrix}
	1_{r-1}+vsu & vs^2\\ -uvu & 1-uvs 
\end{pmatrix}\tp
$$}

This matrix is invertible with $\mu(u,s,v)^{-1}=\mu(u,s,-v)$. If $s\in A^*$, then
$$
\mu(u,s,v)=\begin{pmatrix} 1_{r-1} & 0\\  -u/s & 1 \end{pmatrix}
\begin{pmatrix} 1_{r-1} & vs^2\\  0 & 1 \end{pmatrix}
\begin{pmatrix} 1_{r-1} & 0\\  u/s & 1 \end{pmatrix}
$$ 

The following lemma was proved in \cite{VasBounded}.

{\lem\label{VasersteinMu} (Lemma~2.2 in \cite{VasBounded}) When $r\ge 3$, the matrix $\mu(u,s,v)$ is a product of $7r-3$ elementary transvections in $\GL(r,R)$.}

\smallskip

Now we are ready to prove Proposition~\ref{MakeMonic}.

\smallskip 

Let $S$ be a multiplicative system of lexicographically monic polynomials in $A$. It follows from Lemma~\ref{DimOfLocalisation} and \ref{DimmaxImpliesAFSR} that the ring $A[S^{-1}]$ satisfies $\AFSR_r$, and so does any quotient of $A[S^{-1}]$.

\subsection*{Even case} Take $b\in \IUm_{2r}A$. We need to obtain a lexicographically monic polynomial by $11r-7$ elementary orthogonal transvections. We perform the following steps.
\smallskip

{\bf Step 1.} Make the row $(b_2,\ldots,b_{-1})$ unimodular in $A[S^{-1}]$ by $r-1$ transvections.

\smallskip

Let $\A=\<b_{-r},\ldots,b_{-1}\>\unlhd A[S^{-1}]$. Since $A[S^{-1}]/\A$ satisfies $\AFSR_r$ and the row $(b_1,\ldots,b_r)$ is unimodular in $A[S^{-1}]/\A$, it follows from Lemma~\ref{UseAFSR} that there exist $c_2$,$\ldots$,$c_r\in A$ such that the row $(b_2+c_2b_1,\ldots,b_r+c_rb_1)$ is unimodular in $A[S^{-1}]/\A$. Thus by applying the transvections $T_{i,1}(c_i)$ for $i=2$,$\ldots$,$r$, we make the row $(b_2,\ldots,b_r)$ unimodular in $A[S^{-1}]/\A$ without changing the ideal $\A$. Thus the row $(b_2,\ldots,b_{-1})$ becomes unimodular in $A[S^{-1}]$.

\smallskip

{\bf Step 2.} Make the row $(b_1, b_{-r}\ldots,b_{-1})$ unimodular in $A[S^{-1}]$ by $r-1$ transvections.

\smallskip

Since the row  $(b_2,\ldots,b_{-1})$ unimodular in $A[S^{-1}]$, it follows that the ideal generated by $(b_2,\ldots,b_{-1})$ in $A$ contains a lexicographically monic polynomial. So for some $f_2$,$\ldots$,$f_{-1}\in A$, the polynomial
$$
\sum_{i=2}^{-1} f_ib_i
$$
is lexicographically monic. Multiplying polynomials $f_i$ by a large enough power of $x_1$, we may assume that the polynomial  
$$
b_1+\sum_{i=2}^{-1} f_ib_i
$$
is also lexicographically monic.

Let us now apply the transvections $T_{1,i}(f_i)$ for $i=2$,$\ldots$,$r$. Then the ideal generated by the new $b_1$ and old $b_{-r}$,$\ldots$,$b_{-1}$ contains a lexicographically monic polynomial. However, these transvections do not change the ideal generated by $b_{-r}$,$\ldots$,$b_{-1}$. Hence we actually achieve that the ideal generated by new $b_1$,$b_{-r}$,$\ldots$,$b_{-1}$ contains a lexicographically monic polynomial. Thus the row $(b_1, b_{-r}\ldots,b_{-1})$ becomes unimodular in $A[S^{-1}]$.

\smallskip

{\bf Step 3.} Make the row $(b_1, b_{-r}\ldots,b_{-2})$ unimodular in $A[S^{-1}]$ by $7r-3$ transvections.

\smallskip

Let $\A=\<b_{-r},\ldots,b_{-1}\>\unlhd A[S^{-1}]$. Since $b_1$ is invertible in $A[S^{-1}]/\A$, it follows that there exist $\xi_2$,$\ldots$,$\xi_r\in A[S^{-1}]$ such that $b_i-\xi_ib_1\in \A$ for $i=2,\ldots,r$. Let $s$ be a common denominator of $\xi_i$. So we have $\xi_i=\tfrac{u_i}{s}$. Set 
$$
\tilde{b}_{-1}=b_{-1}+\sum_{i=2}^{r} \xi_ib_{-i}\tp
$$
Since $A[S^{-1}]$ satisfies $\AFSR_r$, it follows from Lemma~\ref{UseAFSR} that there exist $c_{-2}$,$\ldots$,$c_{-r}\in s^2A$ such that every maximal ideal of $A[S^{-1}]$ containing the ideal $\<b_{-r}+c_{-r}\tilde{b}_{-1},\ldots b_{-2}+c_{-2}\tilde{b}_{-1}\>$ contains already the ideal $\<b_{-r},\ldots,b_{-2},\tilde{b}_{-1}\>=\A$. Let $c_i=-s^2v_i$.

Now let $u=(u_r,\ldots,u_2)$, and $v=(v_{-r},\ldots,v_{-2})^t$. We claim that if we multiply the column $b$ by the matrix $H(\mu(u,s,v)^{\t})$, then we make the row $(b_1, b_{-r}\ldots,b_{-2})$ unimodular in $A[S^{-1}]$. Let us prove that.

Assume that some maximal ideal $\M$ of the ring $A[S^{-1}]$ contains the new elements $b_1$,$b_{-r}$,$\ldots$,$b_{-2}$. In terms of the old elements that would be
\begin{gather*}
	(1-uvs)b_1+s^2\left(\sum_{k=2}^r v_{-k}b_k\right)=b_1+s^2\left(\sum_{k=2}^r v_{-k}(b_k-\xi_kb_1)\right)\in\M\tc\\
	b_i-v_i\left(s\left(\sum_{k=2}^r u_k b_{-k}\right)+s^2b_{-1}\right)=b_i+c_i\tilde{b}_{-1}\in\M\qquad i=-r,\ldots,-2\tp
\end{gather*}

By choice of $c_i$ we have $\A\le\M$. Hence $(b_k-\xi_kb_1)\in\M$ for $k=2$,$\ldots$,$r$. Thus $b_1\in\M$. However, by previous step $b_1$ and $\A$ generate a unit ideal. This is a contradiction.

It remains to refer to Lemma~\ref{VasersteinMu}, which implies that $H(\mu(u,s,v)^{\t})\in EO(2r,R)^{\le 7r-3}$. 

\smallskip

{\bf Step 4.} Make the row $(b_1, b_{-r}\ldots,b_{-3})$ unimodular in $A[S^{-1}]$ by $r-1$ transvections.

\smallskip

Since $A[S^{-1}]$ satisfies $\AFSR_r$ and the row $(b_1, b_{-r}\ldots,b_{-2})$ is unimodular in $A[S^{-1}]$, it follows from Lemma~\ref{UseAFSR} that there exist $c_1$,$c_{-r}$,$\ldots$,$c_{-3}\in A$ such that the row $(b_1+c_2b_{-2}, b_{-r}+c_{-r}b_{-2},\ldots,b_{-3}+c_{-3}b_{-2})$ is unimodular in $A[S^{-1}]$. Thus by applying the  transvections $T_{i,-2}(c_i)$ for $i=1$,$-r$,$\ldots$,$-3$, we make the row $(b_1, b_{-r}\ldots,b_{-3})$ unimodular in $A[S^{-1}]$.

\smallskip

{\bf Step 5.} Make $b_{-2}$ lexicographically monic by $r-1$ transvections..

\smallskip

Since the row  $(b_1,b_{-r}\ldots,b_{-3})$ is unimodular in $A[S^{-1}]$, it follows that the ideal generated by $(b_1,b_{-r}\ldots,b_{-3})$ in $A$ contains a lexicographically monic polynomial. So for some $f_1$,$f_{-r}$,$\ldots$,$f_{-3}\in A$, the polynomial
$$
\sum_{i\in\{1,-r,\ldots,-3\}} f_ib_i
$$
is lexicographically monic. Multiplying polynomials $f_i$ by a large enough power of $x_1$, we may assume that the polynomial  
$$
b_{-2}+\sum_{i\in\{1,-r,\ldots,-3\}} f_ib_i
$$
is also lexicographically monic.

Now applying the transvections $T_{-2,i}(f_i)$ for $i=1$,$-r$,$\ldots$,$-3$, we obtain a lexicographically monic polynomial in position $-2$.

\smallskip

\subsection*{Odd case} Take $b\in \IUm_{2r+1}A$. We need to obtain a lexicographically monic polynomial by $12r-8$ elementary orthogonal transvections. The proof is similar to the one given for the even case. However, here we need an additional step at the begining 

\smallskip

{\bf Step 0.} Make the row $(b_1,\ldots,b_r,b_{-r},\ldots,b_{-1})$ unimodular in $A[S^{-1}]$ by $r-1$ transvections.

\smallskip

Let $\A=\<b_1,b_{-r},\ldots,b_{-1}\>\unlhd A[S^{-1}]$. Since $A[S^{-1}]/\A$ satisfies $\AFSR_r$ and the row $(b_2,\ldots,b_r,2b_0)$ is unimodular in $A[S^{-1}]/\A$, it follows from Lemma~\ref{UseAFSR} that there exist $c_2$,$\ldots$,$c_r\in A$ such that the row $(b_2+2c_2b_0,\ldots,b_r+2c_rb_0)$ is unimodular in $A[S^{-1}]/\A$. Thus by applying the transvections $T_{i,0}(c_i)$ for $i=2$,$\ldots$,$r$, we make the row $(b_2,\ldots,b_r)$ unimodular in $A[S^{-1}]/\A$ without changing the ideal $\A$. Thus the row $(b_1,\ldots,b_r,b_{-r},\ldots,b_{-1})$ becomes unimodular in $A[S^{-1}]$.

\smallskip

All the remaining steps in the odd case use only long root transvections and are identical to those in the even case.

{\rem One can notice that the proof above basically repeats the proof of stability theorem for $K_1$-functor given in \cite{SteinStability}.}

\section{Eliminating a variable}
\label{EliminationSection}

In this section, we give the proof of Proposition~\ref{KillTheVariable}. First we need some preparation.

\smallskip

For a commutative ring $R$, we denote by $\Theta(r,R)$ be the following set
$$
\Theta(r,R)=\{M\in M^{r\times r}(R)\colon M^{\t}=-M \text{, and } M \text{ has zeroes at antidiagonal positions}\}\tp
$$

The set $\Theta(r,R)$ forms a group under the addition operation, and the map
$$
\map{U}{\Theta(r,R)}{O(2r,R)}\tc\qquad M\mapsto \begin{pmatrix} 1_r & 0\\ M & 1r\end{pmatrix}
$$
is a group homomorphism.

In the odd case, we denote by $U$ the composition
$$
\xymatrix{
	\Theta(r,R)\ar[r]^U & O(2r,R) \ar@{^{(}->}[r] & O(2r+1,R)\tp
}
$$

{\lem\label{UnipotentRadical} Any matrix from the image of the map $U$ can be expressed as a product of $r(r-1)/2$ elementary orthogonal transvections.}
\begin{proof}
	Suppose that 
	$$
	M=\sum_{1\le j<i\le r} \xi_{i,j}(e_{r+1-i,j}-e_{r+1-j,i})\tp
	$$
	Then
	$$
	U(M)=e+\sum_{1\le i<j\le r} \xi_{i,j}(e_{-i,j}-e_{-j,i})=\prod_{1\le i<j\le r} T_{-i,j}(\xi_{i,j})\tp
	$$
\end{proof}

{\rem From the Chevalley groups point of view, in the even case the image of $U$ is a unipotent radical of a certain parabolic subgroup, and $r(r-1)/2$ is a number of roots in the special part of the corresponding parabolic set of roots.}

{\lem\label{ExistsMUnimodCase} Let $r\ge 2$. Let $b^+=(b_1,\ldots,b_r)^t\in R^r$ be a unimodular column. Let $b^-=(b_{-r},\ldots,b_{-1})^t\in R^r$ be such that $\sum_{i=1}^r b_ib_{-i}=0$. Then there exists a matrix $M\in \Theta(r,R)$ such that $b^-=Mb^+$.}

\begin{proof}
	For a system of linear equations over a ring, existence of a solution is a local property, see, for example, Proposition~1 in \cite{LinEq}. Hence it is enough to consider the case where $R$ is a local ring.
	
	In the local case, since $b^+$ is unimodular, there exists a matrix $g\in \GL(r,R)$ such that $gb^+=e_1$. Set $\tilde{b}^-=(\tilde{b}_{-r},\ldots,\tilde{b}_{-1})^t=(g^{\t})^{-1}b^-$. Note that the first column of $g^{-1}$ is $b^+$. Hence the last row of $(g^{\t})^{-1}$ is transposed $b^+$ written backwards. Hence $\tilde{b}_{-1}=0$.
	
	Further set
	$$
	\tilde{M}=\begin{pmatrix}	
		\tilde{b}_{-r}   &  				& 		 &  \\
		\vdots 		   	 &  			 	& 		 &  \\
		\tilde{b}_{-2} 	 &  			 	& 		 &  \\
		0     			 & -\tilde{b}_{-2}  &\ldots  & -\tilde{b}_{-r}
	\end{pmatrix}\tp
	$$
	Then we have $\tilde{M}e_1=\tilde{b}^-$. It remains to set $M=g^{\t}\tilde{M}g$.
\end{proof}

{\lem\label{ExistsMForskbminus} For any $r\ge 2$, there exists $\tilde{k}\in\N$ such that for any ring $R$, for any columns $b^+=(b_1,\ldots,b_r)^t\in R^r$ and $b^-=(b_{-r},\ldots,b_{-1})^t\in R^r$ such that $\sum_{i=1}^r b_ib_{-i}=0$, and for any $s\in \<b_1,\ldots,b_r\>\unlhd R$, there exists a matrix $M\in \Theta(r,R)$ such that $s^{\tilde{k}}b^-=Mb^+$.}

\begin{proof}
	First we prove that such a number $\tilde{k}$ exists for any given $R$, $b^+$, $b^-$ and $s$. Since the column $b^+$ becomes unimodular in the localisation $R[s^{-1}]$, it follows from Lemma~\ref{ExistsMUnimodCase} that there exists a matrix $M_1\in \Theta(r,R[s^{-1}])$ such that $b^-=M_1b^+$. Let $l_1\in\N$ be such that $M_1=M_2s^{-l_1}$ for some $M_2\in \Theta(r,R)$. Then the equality $s^{l_1} b^-=M_2b^+$ holds in $R[s^{-1}]$. Hence for a large enough $l_2$, we have $s^{l_1+l_2} b^-=s^{l_2}M_2b^+$ in $R$. Set $\tilde{k}=l_1+l_2$ and $M=s^{l_2}M_2$.
	
	Now we show that $\tilde{k}$, actually, can depend only on $r$. Consider the universal situation where
	$$
	R=\Z[b_1,\ldots,b_{-1},a_1,\ldots,a_r]/\<b_1b_{-1}+\ldots+b_rb_{-r}\>\tc\qquad\text{and}\qquad s=a_1b_1+\ldots+a_rb_r\tp
	$$
	Since this case can be specialised to any other, it follows that the number $\tilde{k}$ that works in this case should work in any other case as well.
\end{proof}

{\lem\label{ExistsMIfDivisible} Let $R$ be a commutative Noetherian ring, $s\in R$ and $r\ge 2$. Then there exists $k\in \N$ such that for any column $b^+=(b_1,\ldots,b_r)^t\in R^r$ such that $s\in \<b_1,\ldots,b_r\>$, and for any column $b^-=(b_{-r},\ldots,b_{-1})^t\in R^r$ such that $\sum_{i=1}^r b_ib_{-i}=0$ and all the entries $b_{-i}$ are divisible by $s^k$, there exists a matrix $M\in \Theta(r,R)$ such that $b^-=Mb^+$.}

\begin{proof}
	The annihilators of the elements $s^i$, $l\in\N$ form a ascending chain
	$$
	\Ann s\le \Ann s^2\le \ldots\tp
	$$  
	Since the ring $R$ is Noetherian, it follows that for some $l\in \N$, we have $\Ann s^{l+m}=\Ann s^l$ for any $m\in\N$.
	
	Take $\tilde{k}$ from Lemma~\ref{ExistsMForskbminus}, and set $k=\tilde{k}+l$. Now let  $b^+=(b_1,\ldots,b_r)^t\in R^r$ be such that $s\in \<b_1,\ldots,b_r\>$, and $b^-=(b_{-r},\ldots,b_{-1})^t\in R^r$ such that $\sum_{i=1}^r b_ib_{-i}=0$ and $b^-=s^k (b^-)'$ for some $(b^-)'\in R^k$.
	
	Since $\sum_{i=1}^r b_ib_{-i}=0$, it follows that $\sum_{i=1}^r b_ib_{-i}'\in \Ann s^k=\Ann s^l$. Hence $\sum_{i=1}^r b_i(s^lb_{-i}')=0$. By the choice of $\tilde{k}$ there exists $M\in \Theta(r,R)$ such that $b^-=s^{\tilde{k}}(s^l(b^-)')=Mb^+$. 
\end{proof}

{\lem \label{AddToY} Let $B$ be a commutative Noetherian ring, $A=B[y]$, $r\ge 3$, $b=b(y)=(b_1,\ldots,b_{-1})^t\in \IUm_{2r}A$ resp. $\IUm_{2r+1}A$, $s\in B\cap \<b_1,\ldots,b_{r-1}\>$. Then there exists $m\in\N$ such that
$$
b(y+s^mz)\in EO(2r,A[z])^{\le r(r-1)/2+8r-4}b(y)\;\text{resp.}\; b(y+s^mz)\in EO(2r+1,A[z])^{\le r(r-1)/2+9r-5}b(y)\tp
$$}

\begin{proof}
	First we consider the even case.
	
	Take $k$ from Lemma~\ref{ExistsMIfDivisible} (for $R=A[z]$) and set $m=k+2$. Let $b=(\begin{smallmatrix}
		b^+\\ b^-
	\end{smallmatrix})$, where  $b^+=(b_1,\ldots,b_r)^t\in A^r$, and $b^-=(b_{-r},\ldots,b_{-1})^t\in A^r$. Since $s\in B\cap \<b_1,\ldots,b_{r-1}\>$, by Corollary 2.4 of \cite{VasBounded} there exists a matrix $g(z)\in E(r,A[z])^{\le 8r-4}$ such that $g(z)b^+(y)=b^+(y+s^2z)$. 

	Moreover, it follows from the proof that $g$ is congruent to the identity matrix modulo $z$. Indeed, in Vaserstein's proof he reduces the column $b(y+s^2z)$ to $b(y)$ in two steps. At the first step he makes the first $r-1$ components of the column equal to those of $b(y)$ using matrix $\mu(u,s,v)$ (see Definition \ref{mu} of the present paper) with all the entries of the column $v$ being divisible by $z$; hence this matrix is congruent to the identity matrix modulo $z$. It follows that the difference between $b_r(y)$ and the last entry of the column obtained on the first step is divisible by $z$. At the second step Vaserstein uses addition operations to make the last entry equal to $b_r(y)$. In a view of the above, the coefficients of these additions must be divisible by $z$. Thus the matrix of the whole transformation is congruent to the identity matrix modulo $z$.
	
	Therefore, for some $\tilde{b}^-\in A[z]^r$ we have
	$$
	b'=\begin{pmatrix} b^+(y+s^mz) \\ \tilde{b}^-\end{pmatrix}=H(g(s^kz))b\in EO(2r,A[z])^{\le 8r-4}b\tp
	$$
	In addition, $\tilde{b}^-$ is congruent to $b^-(y)$ modulo $s^k$. Since both $b(y+s^mz)$ and $b'$ are isotropic, and $b^-(y+s^mz)-\tilde{b}^-$ is divisible by $s^k$, it follows by Lemma~\ref{ExistsMIfDivisible} that there exists a matrix $M\in \Theta(r,A[z])$ such that $b^-(y+s^mz)-\tilde{b}^-=Mb^+(y+s^mz)$. Here we used that $s\in B$, so the shift of the variable does not change the fact that  $s\in \<b_1,\ldots,b_{r-1}\>$. Hence by Lemma~\ref{UnipotentRadical} we have
	$$
	b(y+s^mz)=U(M)b'\in EO(2r,A[z])^{\le r(r-1)/2+8r-4}b
	$$
	
	In the odd case, where we have 
	$$
	b=\begin{pmatrix} b^+ \\ b_0 \\ b^-\end{pmatrix}\tc
	$$
	the proof is the similar. However, after creating $b^+(y+s^mz)$ at the top and before creating $b^-(y+s^mz)$ at the bottom, we have to create $b_0(y+s^mz)$ in the middle by the transvections $T_{-i,0}(\xi_i)$, $i=1,\ldots,r-1$ with suitable $\xi_i$. Here we again use that $s\in \<b_1,\ldots,b_{r-1}\>$.
\end{proof}

{\lem \label{PrepareToAddToY} Let $B$ be a commutative ring, $P_1$,$\ldots$,$P_m$ be distinct maximal ideals in $B$, $A=B[y]$, $r\ge 3$, $b=b(y)=(b_1,\ldots,b_{-1})^t\in \IUm_{2r}A$ resp. $\IUm_{2r+1}A$ such that $b_{-2}$ is monic. Then there exists a column
$$
b^{(1)}\in EO(2r,A)^{\le \ceil{{r-1\over 2}}+2}b\qquad\text{resp.}\qquad EO(2r+1,A)^{\le  \ceil{{r-1\over 2}}+3}b
$$
such that $b^{(1)}_{-2}$ is monic and
$$
\left(\<b^{(1)}_1,\ldots b^{(1)}_{r-1}\>\cap B\right)\sm \bigcup_{i=1}^m P_i\ne \emp\tp
$$
}
\begin{proof}
		Set 
		$$
		R=B/\left(\bigcap_{i=1}^m P_i\right)=\prod_{i=1}^m B/P_i\tp
		$$
		
		First we show that the last condition on $b^{(1)}$ holds if $b^{(1)}_1$ monic and the row $(b^{(1)}_1,\ldots,b^{(1)}_{r-1})$ becomes unimodular in $R[y]$.
		
		Let $c_1$,$\ldots$,$c_{r-1}\in A$ be such that $c_1b^{(1)}_1+\ldots+c_{r-1}b^{(1)}_{r-1}\equiv 1\mod P_i$ for every $i$.
		
		Set $f=c_2b^{(1)}_2+\ldots+c_{r-1}b^{(1)}_{r-1}$. Then $b^{(1)}_1$ and $f$ are coprime in $B/P_i[y]$ for every $i$.
		
		Since $b^{(1)}_1$ is monic, it follows that the resultant $\res(b^{(1)}_1,f)$ modulo $P_i$ is equal to the resultant of $b^{(1)}$ taken modulo $P_i$ and $f$ taken modulo $P_i$ (even if $f$ modulo $P_i$ has smaller degree).
		
		Therefore, we have 
		$$
		\res(b^{(1)}_1,f)\in \left(\<b^{(1)}_1,\ldots b^{(1)}_{r-1}\>\cap B\right)\sm \bigcup_{i=1}^m P_i\tp
		$$
		
		Thus it remains to prove that a given column $b\in \IUm_{2r}A$ resp. $\IUm_{2r+1}A$, with $b_{-2}$ being monic, can be transformed by $\ceil{{r-1\over 2}}+2$ resp. $\ceil{{r-1\over 2}}+3$ transvections so that $b_1$ becomes monic, $b_{-2}$ remains monic, and the row $(b_1,\ldots,b_{r-1})$ becomes unimodular in $R[y]$.
		
		\smallskip
		
		First we consider the even case. Here we perform the following steps.
		
		\smallskip
		
		{\bf Step 1.} Make the polynomial $b_1$ monic and the row $(b_1,\ldots,b_{-2})$ unimodular in $R[y]$ by the transvection $T_{1,-2}(\xi)$.
		
		\smallskip
		
		Since $R$ is a product of fields and $b_{-2}$ is monic, it follows that the ring $R[y]/\<b_1,b_3,\ldots,b_{-2}\>$ is semilocal; hence this ring satisfies $\SR_2$. Then there exists $\tilde{\xi}$ such that the row $(b_1,b_2+\tilde{\xi}b_{-1},b_3,\ldots,b_{-2})$ is unimodular in $R[y]$. Therefore, if we take 
		$$
		\xi=\tilde{\xi}+y^Kb_{-2}
		$$
		for some $K\in\N$, then we guarantee that the row $(b_1,\ldots,b_{-2})$ becomes unimodular in $R[y]$. It remains to notice that if $K$ is large enough, then we also make $b_1$ monic.
		
		\smallskip
		
		{\bf Step 2.} If $r$ is even, make the row $(b_1,\ldots,b_{-3})$ unimodular in $R[y]$ by the transvection $T_{-3,-2}(\xi)$. If $r$ is odd, then skip this step.
		
		\smallskip
		
		Since $b_1$ is now monic, it follows that the ring $R[y]/\<b_1,\ldots,b_{-4}\>$ is semilocal; hence this ring satisfies $\SR_2$. Therefore, there exists an element $\xi$ with the required property. 
		
		\smallskip
		
		{\bf Step 3.} Make the row $(b_1,\ldots,b_r)$ unimodular in $R[y]$ by $\floor{{r-1\over 2}}$ transvections.
		
		\smallskip
		
		It is enough to show that if the row $(b_1,\ldots,b_{-i+2})$ is unimodular in $R[y]$, then we can make the row  $(b_1,\ldots,b_{-i})$ unimodular in $R[y]$ by the transvection $T_{i-1,-i+2}(\xi)$, where $4\le i\le r+1$. Here if $i=r+1$, then by  $(b_1,\ldots,b_{-i})$ we mean $(b_1,\ldots,b_r)$.
		
		Since $b_1$ is now monic, it follows that the ring $R[y]/\<b_1,\ldots,b_{i-3},b_{i}\ldots,b_{-i}\>$ is semilocal. Let $\M_1$,$\ldots$, $\M_l$ be its maximal ideals. Then it is enough to choose $\xi$ to be contained in those and only those $\M_j$ that do not contain the image of the ideal $\<b_{i-2},b_{i-1}\>$. 
		
		\smallskip
		
		{\bf Step 4.} Make the row $(b_1,\ldots,b_{r-1})$ unimodular in $R[y]$ by the transvection $T_{r-1,r}(\xi)$.
		
		\smallskip
		
		This step is similar to Step 2. 
		 
		\smallskip
		
		It remains to note that neither of those steps change $b_{-2}$. Hence it remains unimodular.
		
		\smallskip
		
		The proof in the odd case is similar. However, here we need an additional step at the begining 
		
		\smallskip
		
		{\bf Step 0.} Make the row $(b_1,\ldots,b_r,b_{-r},\ldots,b_{-1})$ unimodular in $R[y]$ by the transvection $T_{1,0}(\xi)$.
		
		\smallskip
		
		This step is similar to Steps 2 and 4. 
\end{proof}

Now we are ready to prove Proposition~\ref{KillTheVariable}.

\smallskip

For simplicity, we will write
$$
\xymatrix{
v\ar[r]^{N_1}_{N_2} & w
}
$$
instead of
$$
w\in EO(2r,R)^{\le N_1}v\qquad\text{resp.}\qquad w\in EO(2r+1,R)^{\le N_2}v\tc
$$
where $v$ and $w$ are columns.

\smallskip

Applying Lemmas~\ref{BSInduction} and~\ref{PrepareToAddToY} $\delta$ times, we obtain elements $s_1$,$\ldots$,$s_{\delta}\in B$ and columns $b=b^{(0)}$,$ b^{(1)}$,$\ldots$,$b^{(\delta)}\in \IUm_{2r} A$ resp. $\IUm_{2r+1} A$ such that, firstly, 
$$
\xymatrix{
	b^{(i)}\ar[rr]^{\ceil{{r-1\over 2}}+2}_{\ceil{{r-1\over 2}}+3} & & b^{(i+1)} & &  i=0,\ldots,\delta-1\tc
}
$$
secondly, $s_i\in\<b^{(i)}_1,\ldots, b^{(i)}_{r-1}\>$ for $i=1,\ldots,\delta$, and thirdly, $\BSdim B/(s_1,\ldots s_{i+1})<\BSdim B/(s_1,\ldots,s_i)$ for $i=0,\ldots,\delta-1$. In particular, the elements $s_1$,$\ldots$,$s_{\delta}$ generate the unit ideal.

By Lemma~\ref{AddToY} we have
$$
\xymatrix{
	b^{(i)}(y)\ar[rrr]^{r(r-1)/2+8r-4}_{r(r-1)/2+9r-5} & & & b^{(i)}(y+s_i^{m_i}z)
}
$$
in $A[z]$.

Therefore, we have the following chain of transformations in $A[z_1,\ldots,z_{\delta}]$:

\begin{multline*}
\xymatrix{
	b=b^{(0)}(y)\ar[r] & b^{(1)}(y)\ar[r] & b^{(1)}(y+s_1^{m_1}z_1)\ar[r] & b^{(2)}(y+s_1^{m_1}z_1)\ar[r] & \ldots}\\ \xymatrix{b^{(\delta)}(y+s_1^{m_1}z_1+\ldots+s_{\delta-1}^{m_{\delta-1}}z_{\delta-1})\ar[r] & b^{(\delta)}(y+s_1^{m_1}z_1+\ldots+s_\delta^{m_\delta}z_\delta)\tp
}
\end{multline*}

Thus we have
$$
\xymatrix{
	b(y)\ar[rrrrr]^{\!\!\!\!\!\!\!\!\!\!\!\!\!\!\!\!\!\!\!\!\!\delta\left(\tfrac{r(r-1)}{2}+\ceil{\tfrac{r-1}{2}}+8r-2\right)}_{\!\!\!\!\!\!\!\!\!\!\!\!\!\!\!\!\!\!\!\!\!\delta\left(\tfrac{r(r-1)}{2}+\ceil{\tfrac{r-1}{2}}+9r-2\right)} & & & & & b^{(\delta)}(y+s_1^{m_1}z_1+\ldots+s_\delta^{m_\delta}z_\delta)\tp
}
$$

Since the elements $s_1$,$\ldots$,$s_{\delta}$ generate the unit ideal, it follows that so do the elements $s_1^{m_1}$,$\ldots$,$s_{\delta}^{m_{\delta}}$. Specializing the indeterminates $z_i$ to elements in $yB$, we make $y+s_1^{m_1}z_1+\ldots+s_\delta^{m_\delta}z_\delta$ equal to zero; this concludes the proof of Proposition~\ref{KillTheVariable}.

\section*{Acknowledgment}

The author is a participant of a scientific group that won "Leader" grant by "BASIS" foundation in 2020, grant \#20-7-1-27-4. Research is also supported by «Native towns», a social investment program of PJSC ''Gazprom Neft''.


\end{document}